\documentclass[letterpaper, 10 pt, conference]{ieeeconf}
\IEEEoverridecommandlockouts
\makeatletter

\let\proof\@undefined
\let\endproof\@undefined

\makeatother

\usepackage{graphicx}
\usepackage{bbm}
\usepackage{dsfont}
\usepackage{amsthm}
\usepackage{amsmath}
\usepackage{amssymb}
\usepackage[mathscr]{eucal}
\usepackage{url}
\usepackage{algorithm}
\usepackage{algorithmic}
\usepackage[font=scriptsize,caption=false]{subfig}
\usepackage{color}

\newcommand{\cf}[1]{(cf.~#1)}
\newcommand{\set}[1]{\left\{ #1 \right\}}
\newcommand{\paren}[1]{\left( #1 \right)}

\newcommand{\vf}[1]{\frac{\partial\ }{\partial #1}}

\newcommand{\bd}{\partial}

\newcommand{\rank}[1]{\text{rank}\,#1}
\newcommand{\td}[1]{\tilde{#1}}

\newcommand{\eqn}[1]{\begin{equation*} #1 \end{equation*}}
\newcommand{\eqnn}[1]{\begin{equation} #1 \end{equation}}
\newcommand{\eqnal}[1]{\begin{equation*}\begin{aligned}[b] #1 \end{aligned}\end{equation*}}

\newcommand{\sm}{\setminus}
\newcommand{\into}{\rightarrow}

\newcommand{\inc}{\hookrightarrow}

\newcommand{\R}{\mathbb{R}}

\newcommand{\N}{\mathbb{N}}

\newcommand{\e}{\mathscr}
\newcommand{\eps}{\epsilon}
\newcommand{\vphi}{\varphi}

\newcommand{\rk}{\operatorname{rank}}
\newcommand{\Int}[1]{\operatorname{Int}#1}

\newcommand{\Pmap}{Poincar\'{e} map}

\newtheorem{theorem}{Theorem}
\newtheorem{proposition}{Proposition}
\newtheorem{corollary}{Corollary}
\newtheorem{lemma}{Lemma}
\newtheorem{assumption}{Assumption}

\newtheorem{remark}{Remark}
\newtheorem{example}{Example}
\newtheorem{definition}{Definition}

\newcommand{\defn}[1]{\begin{definition} #1 \end{definition}}

\newcommand{\ex}[1]{\begin{example} #1 \end{example}}
\newcommand{\rem}[1]{\begin{remark} #1 \end{remark}}
\newcommand{\assump}[1]{\begin{assumption} #1 \end{assumption}}
\newcommand{\thm}[1]{\begin{theorem} #1 \end{theorem}}
\newcommand{\lem}[1]{\begin{lemma} #1 \end{lemma}}
\newcommand{\cor}[1]{\begin{corollary} #1 \end{corollary}}
\newcommand{\pf}[1]{\begin{proof} #1 \end{proof}}
\newcommand{\prop}[1]{\begin{proposition} #1 \end{proposition}}

\begin{document}

\title{Dimension Reduction Near Periodic Orbits of Hybrid Systems}

\author{Samuel Burden, Shai Revzen and S. Shankar Sastry%
\thanks{S. Burden and S. Sastry are with the Department of Electrical Engineering and Computer Sciences,
        University of California at Berkeley, CA, USA
        {\tt\small sburden,sastry@eecs.berkeley.edu}}%
\thanks{S. Revzen is with the Department of Electrical and Systems Engineering, University of Pennsylvania,
        Philadelphia, PA, USA
        {\tt\small shrevzen@seas.upenn.edu}}%
}

\maketitle

\begin{abstract}
When the Poincar\'{e} map associated with a periodic orbit of a hybrid dynamical system has constant-rank iterates, we demonstrate the existence of
a constant-dimensional invariant subsystem near the orbit 
which attracts all nearby trajectories in finite time.
This result shows that the long-term behavior of a hybrid model with a large number of degrees-of-freedom may be governed by a low-dimensional smooth dynamical system.
The appearance of such simplified models 
enables the translation of analytical tools from smooth systems---such as Floquet theory---to the hybrid setting
and provides a bridge between the efforts of biologists and engineers studying legged locomotion.
\end{abstract}

\section{Introduction}

Dynamic multi-legged locomotion presents a daunting control task.
A large number of degrees-of-freedom (DOF) must be rapidly and precisely coordinated in the face of state and environmental uncertainty.  
The ability of individual limbs to exert forces on the body varies intermittently with ground contact, body posture, and the efforts of other appendages.
Finally, the motion itself affects sensor measurements, complicating pose estimation.
In spite of these difficulties, animals at all levels of complexity have mastered the art of rapid legged locomotion over complex terrain at speeds far exceeding those of comparable robotic platforms~\cite{DickinsonFarley2000, Raibert1986, SaranliBuehler2001, KimClark2006}.

\begin{figure}
\centering
\includegraphics[scale=0.95]{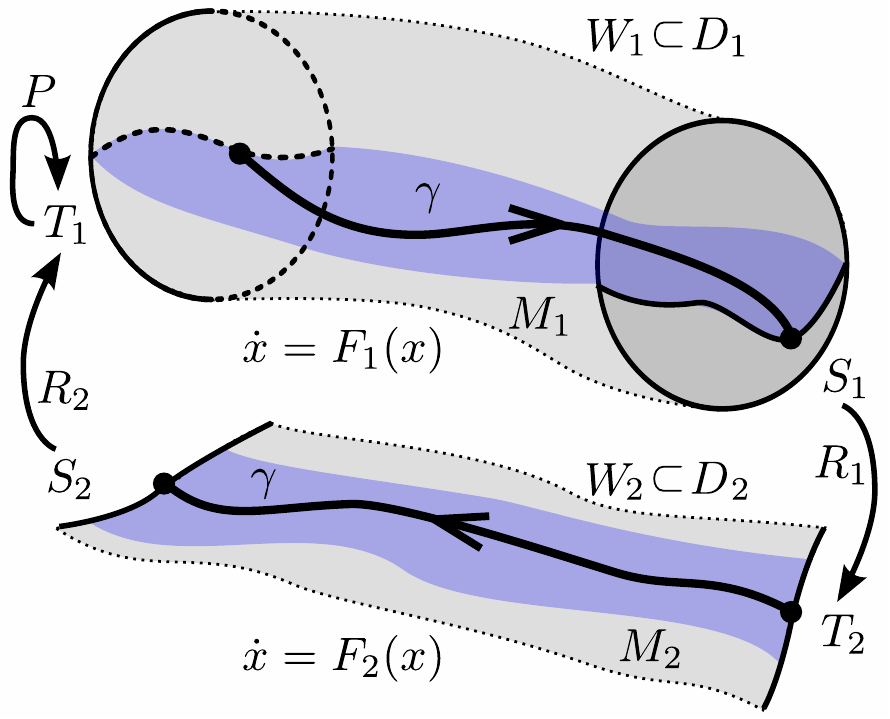}
\caption{
Illustration of the main result in a two-domain system:  
Whenever iterates of the \Pmap\, $P$ associated with a periodic orbit $\gamma$ of a hybrid dynamical system are constant-rank,
there is a constant-dimensional smooth submanifold $M_j$ in each domain $D_j$ which is invariant under the flow of the vector field $F_j$ and which attracts all trajectories starting in an open set $W_j\subset D_j$ containing $\gamma$ in finite time.
}
\label{fig:result}
\end{figure}

Numerous architectures have been proposed to explain how animals control their limbs.
For steady-state locomotion, most posit a principle of coordination, synergy, symmetry or synchronization, and there is a surfeit of neurophysiological data to support these hypotheses~\cite{Grillner1985, CohenHolmesRand1982, GolubitskyEtAl1999, TingMacph2005}.
In effect, the large number of DOF available to an animal are collapsed during regular motion to a low-dimensional dynamical attractor that may be captured by a \emph{template} model embedded within a higher-dimensional model \emph{anchored} to the animal's morphology~\cite{HolmesEtAl2006, FullKoditschek1999}.
In this view, only a few parameters like frequency and coupling strength are required to describe the dynamics of any particular periodic gait over a broad range of animal morphologies, offering a tantalizing target for experimental biologists. 
Were the dynamics of legged animals smooth as a function of position and momentum, Floquet theory~\cite{floquet1883edl} provides a canonical form for the structure of the stability basin of a limit cycle~\cite{Guckenheimer1975, GuckenheimerHolmes1983}.
In such a canonical form, the template may appear as an invariant attractor of the linearized dynamics and be amenable to quantitative measurement~\cite{Revzen2009, RevzenGuckenheimer2011}.
A substantial motivation for the present work has been to provide a theoretical framework for applying this empirical approach to study legged locomotion.
The dynamics of legged locomotion are rarely smooth due to intermittent contact of limbs with the substrate, so we have generalized this approach to be aplicable to a class of non-smooth systems called \emph{hybrid dynamical systems}.

We relegate a formal definition of the class of hybrid systems under consideration to Section~\ref{sec:HDS}. 
Informally, hybrid dynamical systems are comprised of differential equations written over disparate domains together with rules for switching between the domains.
Of particular interest to us are periodic orbits of such systems.
From a modeling viewpoint, a stable hybrid periodic orbit provides a natural abstraction for the dynamics of steady-state legged locomotion.  
This approach has been widely adopted, generating a variety of models of bipedal~\cite{McGeer1990, GrizzleAbba2002, CollinsRuina2005} and multi-legged~\cite{GhigliazzaAltendorfer2003, ProctorHolmes2008} locomotion as well as some general control-theoretic techniques for composition~\cite{KlavinsKoditschek2002}, coordination~\cite{HaynesCohen2009}, and stabilization~\cite{WesterveltGrizzle2003, MorrisGrizzle2009, ShiriaevFreidovich2010} of such models.
In certain cases, it has been possible to formally embed a low-dimensional abstraction in a higher-dimensional physically-realistic model~\cite{PoulakakisGrizzle2009, AmesGregg2006}.

This paper provides a conceptual link between formal analysis of hybrid periodic orbits and the dramatic dimension reduction observed empirically in successful legged locomotors.  
Under the condition that iterates of the \Pmap\, associated with a periodic orbit are constant rank, we demonstrate the existence of a constant-dimensional invariant subsystem which attracts all nearby trajectories in finite time.  
Analogous results for smooth dynamical systems typically impose stringent assumptions on the dynamics such as exact symmetries \cf{\S 8.9 in~\cite{MarsdenRatiu1999}} or timescale separation \cf{Chapter 4 in~\cite{GuckenheimerHolmes1983}}.
In contrast, the results of this paper imply that hybrid dynamical systems may exhibit dimension reduction near periodic orbits solely due to the interaction of the switching dynamics with the smooth flow.

\subsection*{Organization}

The hybrid systems we consider are constructed using switching maps defined between boundaries of smooth dynamical systems.  
The behavior of such systems can be studied by 
alternately applying flows and maps. 
Thus, we begin in Section~\ref{sec:SDS} by developing several results which provide canonical forms for the behavior of flows and maps near periodic orbits and fixed points, respectively.  
Then, we define hybrid systems in Section~\ref{sec:HDS} and use these results to characterize the dynamics near their periodic orbits.
Examples are presented in Section~\ref{sec:ex} and implications of the results for the design and analysis of legged locomotors are explored in Section~\ref{sec:disc}.

\section{Smooth Dynamical Systems}\label{sec:SDS}

This section contains three technical results used in the proof of the Theorem of Section~\ref{sec:HDS}.  
The first two results concern smooth dynamical systems\footnote{For notational convenience, we work with objects which possess continuous derivatives of all orders.  However, the results in this paper are valid if we only assume continuous differentiability.} and may be found in textbooks, hence we state them without proof.
The third establishes, under a non-degeneracy condition, a canonical form for the invariant set of a smooth map near a fixed point. 
A reader interested in the main result of this paper may proceed to Section~\ref{sec:HDS} and refer to this section as needed.

\subsection{Differential Geometry}

We assume familiarity with the tools and terminology of differential geometry.
If any of the concepts we discuss are unfamiliar, we refer the reader to~\cite{MarsdenRatiu1999, Lee2002} for more details.

A \emph{smooth $n$-dimensional manifold} $M$ \emph{with boundary} $\bd M$ is an $n$-dimensional topological manifold covered by a collection of \emph{smooth coordinate charts} $\set{(U_\alpha,\vphi_\alpha)}_\alpha$ where $U_\alpha\subset M$ is open and $\vphi_\alpha:U_\alpha\into H^n$ is a homeomorphism where $H^n :=\set{(y_1,\dots,y_n)\in\R^n : y_n \ge 0}$ is the upper half-space.
The charts are \emph{smooth} in the sense that $\vphi_\alpha\circ\vphi_\beta^{-1}$ is a diffeomorphism over $\vphi_\beta(U_\alpha\cap U_\beta)$ for all pairs $\alpha,\beta$.
The boundary $\bd M\subset M$ contains those points which are mapped to the plane $\set{(y_1,\dots,y_n)\in\R^n : y_n = 0}$ in some chart.
We say $S\subset M$ is a \emph{smooth embedded $k$-dimensional submanifold} if near every $x\in S$ there is a smooth coordinate chart $(U_x,\vphi_x)$ so that 
\eqn{\vphi_x(S\cap U_x) \subset \set{y \in \R^n : y_{k+1} = \cdots = y_n = 0}.}
These charts yield \emph{slice coordinates} for the submanifold, and the integer $n-k$ is the \emph{codimension} of $S$. 
It is a straightforward consequence that $\bd M$ is a smooth embedded submanifold without boundary and has codimension~1.  We denote the interior of $M$ by $\Int(M) := M\sm\bd M$.

Each $x\in M$ has an associated \emph{tangent space} $T_x M$, and the disjoint union of the tangent spaces at each point is the \emph{tangent bundle} $TM := \coprod_{x\in M}T_x M$; note that any element in $TM$ may be regarded as a pair $(x,v)$ where $x\in M$ and $v\in T_x M$.
We let $\e{T}(M)$ denote the space of \emph{smooth vector fields} on $M$, i.e. smooth maps $G:M\into TM$ for which $G(x) = (x,v)$ for some $v\in T_x M$ and all $x\in M$.
It is a fundamental result that any $G\in\e{T}(M)$ determines an ordinary differential equation on the manifold which may be solved globally to obtain a \emph{maximal flow} $\phi:\e{F}\into M$ where $\e{F}\subset \R\times M$ is the \emph{maximal flow domain} \cf{Theorem 17.8 in~\cite{Lee2002}}.  
This flow has several important properties which we will use repeatedly;
let $\e{F}^x:=\set{t\in\R : (t,x)\in\e{F}}$ and 
$\e{F}_t:=\set{x\in M : (t,x)\in\e{F}}$.
First, for any initial condition $x\in M$, $\phi(\cdot,x):\e{F}^x\into M$ is the \emph{maximal integral curve} of $G$ passing through $x$, i.e. $\vf{t}\phi(t,x) = G(\phi(t,x))$ for all $t\in\e{F}^x$; we will alternately refer to integral curves as \emph{trajectories}.
Second, for any smooth embedded submanifold $S\subset M$ and $t\in\R$ for which $\set{t}\times S\subset\e{F}$, $\set{\phi(t,x) : x\in S}\subset M$ is an embedded submanifold that is diffeomorphic to $S$.

If $f:M\into N$ is a smooth map between smooth manifolds,
then at each $x\in M$ there is an associated linear map $f_*(x):T_xM\into T_{f(x)}N$ called the \emph{pushforward}.
Globally, the pushforward is a smooth map $f_* : TM\into TN$.
In coordinates, it is the familiar Jacobian matrix.
The \emph{rank} of a smooth map $f:M\into N$ at a point $x\in M$ is defined $\rk_x f := \rk f_*(x)$.  
If $\rk_x f = r$ for all $x\in M$, we simply write $\rk f = r$.
If $\rk f = \dim M$ and $f$ is a homeomorphism onto its image, then $f$ is a \emph{smooth embedding}, and the image of $f$ is a smooth embedded submanifold.
In this case, any smooth vector field $G\in\e{T}(M)$ may be pushed forward to a unique smooth vector field $f_* G\in\e{T}(f(M))$.
A vector field $G\in\e{T}(M)$ is \emph{transverse} to a $k$-dimensional embedded submanifold $S$ at $x\in S$ if, in slice coordinates $(U_x,\vphi_x)$ near $x$, the $(k+j)^{th}$ coordinate of $\vphi_*G$ is non-zero for some $j$ between 1 and $n-k$; otherwise $G$ is \emph{tangent}.
If $S\subset\bd M$, $G$ is \emph{inward-pointing} if the $n^{th}$ coordinate of $\vphi_*G$ is positive and \emph{outward-pointing} if it is negative.

With these preliminaries established, we are in a position to define one of the main dynamical objects of interest in this paper.

\defn{A \emph{smooth dynamical system} is a pair $(M,G)$: 
\begin{itemize}
\item[$M$] is a smooth manifold with boundary $\bd M$;
\item[$G$] is a smooth vector field on $M$, i.e. $G\in\e{T}(M)$.
\end{itemize}}

\subsection{Flows Between Surfaces}

We review the fact that the flow near a trajectory passing transversally between two surfaces has a simple form \cf{Chapter 11.2 in~\cite{HirschSmale1974}}.
In particular, nearby trajectories can be obtained from an embedding of a product manifold.
This will be the prototype for the dynamics near a periodic orbit in one domain of a hybrid system.  

\lem{\label{lem:sec}
Let $(M,G)$ be a smooth dynamical system and $\phi:\e{F}\into M$ its maximal flow.
Suppose $T,S\subset M$ are smooth embedded submanifolds, $S$ has codimension~1, $\phi(\alpha,\xi)\in S$ for some $\alpha>0$ and $\xi\in T$ where $G$ is transverse to $T$ at $\xi$ and $S$ at $\phi(\alpha,\xi)$.  Then we have the following consequences:
\begin{enumerate}
  \item[(i)] there is a neighborhood $U\subset M$ containing $\xi$ and a smooth map $\eta:U\into\R$ so that $\eta(\xi) = \alpha$ and for all $x\in U$, $\eta(x) > 0$ and $\phi(\eta(x),x)\in S$;  $\eta$ is called the \emph{time-to-impact} map;
  \item[(ii)] with $V := U\cap T$, the map $\psi:[0,1]\times V\into M$ with 
    \eqn{\psi(\sigma,v) := \phi(\eta(v)\sigma,v)}
    is a smooth embedding into $M$ whose image contains the trajectory $\gamma = \set{\phi(t,\xi) : 0\le t \le \alpha}$.
\end{enumerate}
}

\rem{This lemma is applicable when $T,S\subset\bd M$, which will be relevant in the study of hybrid systems.}

\subsection{Gluing Flows}\label{sec:glue}

In this section, we provide a method for gluing two smooth dynamical systems together along their boundaries to obtain a new smooth system;  
this construction uses basic results from differential topology \cf{Theorem 8.2.1 in~\cite{Hirsch1976}}.
We will use this construction in Section~\ref{sec:HDS} to attach distinct hybrid domains to one another.  

\lem{\label{lem:glue}
Suppose $(M_1,G_1),(M_2,G_2)$ are smooth $n$-dimensional dynamical systems, $\vphi:\bd M_1\into \bd M_2$ is a diffeomorphism, $G_1$ is outward-pointing along $\bd M_1$ and $G_2$ is inward-pointing along $\bd M_2$.  
Then the topological quotient $M := \frac{M_1\coprod M_2}{\bd M_1 \simeq\bd M_2}$ can be made into a smooth manifold for which (i) the inclusions $M_j\inc M$ are smooth embeddings and (ii) there is a smooth vector field $G\in\e{T}(M)$ that restricts to $G_j$ on $M_j$, $j = 1,2$.
}

\rem{The smooth structure described in Lemma~\ref{lem:glue} is unique.
Further, if $\td{G}_1$ and $\td{G}_2$ are any other smooth vector fields on $M_1$ and $M_2$ which satisfy the hypotheses of the Lemma, the corresponding quotient $\td{M}$ is diffeomorphic to $M$ \cf{Chaper 8 in~\cite{Hirsch1976}}.
}

\subsection{Invariant Set of a Smooth Map Near a Fixed Point}\label{sec:maps}

In studying hybrid dynamical systems, we encounter smooth maps $f:M\into M$ which are not diffeomorphisms.
Viewing iteration of $f$ as a discrete dynamical system, we wish to study the behavior of these iterates near a fixed point $f(\xi) = \xi$.
Note that if $f$ has constant rank equal to $k\in\N$, then its image $f(M)\subset M$ is an embedded $k$-dimensional submanifold near $\xi$ by the Rank Theorem \cf{Theorem 7.13 in~\cite{Lee2002}}.  
With an eye toward dimension reduction, one might hope that the composition $(f\circ f):M\into f(M)$ 
is also constant-rank, but this is not generally true\footnote{Consider the map $f:\R^2\into\R^2$ defined by $f(x,y) := (x^2,x)$.}.
If it is true that iterates of $f$ are eventually constant-rank near the fixed point $\xi$, then one can study the behavior of these iterates by restricting the domain to a lower-dimensional submanifold.

\lem{\label{lem:nondeg}
Let $M$ be a smooth manifold, $f:M\into M$ a smooth map with $f(\xi) = \xi$ for some $\xi\in M$, suppose the rank of $f$ is bounded above by $n\in\N$, and suppose the composition of $f$ with itself $n$ times, $f^n$, has constant rank equal to $r\in\N$ on a neighborhood of $\xi$.
Then $f^n(M)$ is an $r$-dimensional embedded submanifold near $\xi$ and there are neighborhoods  $U,V\subset f^n(M)$ containing $\xi$ for which $f$ maps $U$ diffeomorphically onto $V$.
}

In the proof of Lemma~\ref{lem:nondeg}, we make use of an elementary fact from linear algebra.
The result is easily obtained by passing to the Jordan form.
\prop{\label{prop:la1}
If $A\in\R^{m\times m}$ and $\rk A \le n$, then $\rk(A^{2n}) = \rk(A^n)$.
}

\pf{{\it (of Lemma~\ref{lem:nondeg})}
By the Rank Theorem \cf{Theorem 7.13 in~\cite{Lee2002}}, there is a neighborhood $N\subset M$ of $\xi$ for which $\Sigma := f^n(N)$ is an $r$-dimensional embedded submanifold
and by Proposition~\ref{prop:la1} we have
\eqnal{
  \rk\paren{f^n|_\Sigma}_*(\xi) & = \rk (f^n\circ f^n)_*(\xi) 
  \\
  & = \rk f^n_*(\xi).
}
Therefore $\paren{f^n|_\Sigma}_*:T_\xi\Sigma\into T_\xi\Sigma$ is a bijection, so by the Inverse Function Theorem \cf{Theorem 7.10 in~\cite{Lee2002}}, there is a neighborhood $W\subset \Sigma$ containing $\xi$ so that $f^n(W)\subset\Sigma$ and $f^n|_{W}:W\into f^n(W)$ is a diffeomorphism.

By continuity of $f$, there is a neighborhood $L\subset N$ containing $\xi$ for which $f(L)\subset N$ and $f^n(L)\subset W$. 
The set 
$U := f^n(L)$
is a neighborhood of $\xi$ in $\Sigma$.
Further, we have 
\eqnal{
  f(U) & = f\circ f^n(L) = f^n\circ f(L)\subset\Sigma.
}
The restriction $f^n|_U:U\into f^n(U)$ is a diffeomorphism since $U \subset W$, whence $f|_U$ is a diffeomorphism onto its image, $V := f(U)\subset\Sigma$.
}

\section{Hybrid Dynamical Systems}\label{sec:HDS}

We describe a class of hybrid systems useful for modeling legged locomotion, then restrict our attention to the behavior of such systems near periodic orbits.  
It was shown in~\cite{WendelAmes2010} that the Poincar\'{e} map of a hybrid system is generally not full rank.  
We explore the geometric consequences of this rank loss and demonstrate,
under a non-degeneracy condition, the existence of a smooth invariant subsystem which attracts all nearby trajectories in finite time.  

\subsection{Hybrid Differential Geometry}

For our purposes, it is expedient to define hybrid dynamical systems over disjoint unions of smooth manifolds.  

\defn{\label{def:hman}
A \emph{smooth hybrid manifold} is a finite disjoint union of connected smooth manifolds $M = \coprod_{j\in J} M_j$.}

\rem{The dimensions of the constituent manifolds are not required to be equal.}

Differential geometric constructions
which are confined to a single manifold have natural generalizations to such spaces, and we will prepend the modifier ``hybrid'' to make it clear when this generalization is being invoked.
For instance, the \emph{hybrid tangent bundle} $TM$ is the disjoint union of the tangent bundles $TM_j$,
the \emph{hybrid boundary} $\bd M$ is the disjoint union of the boundaries $\bd M_j$, and 
a \emph{hybrid open set} $U\subset M$ is obtained from a disjoint union of open sets $U_j\subset M_j$.
Generalizing maps between manifolds requires more care, hence we provide explicit definitions.

\assump{To simplify the exposition, we henceforth assume all manifolds and maps between manifolds are smooth.}

\defn{\label{def:hmap} 
A \emph{hybrid map} 
\eqn{f:\coprod_{j\in J} M_j\into\coprod_{\ell\in L} N_\ell}
between hybrid manifolds restricts to a map $f|_{M_j} : M_j\into N_\ell$, some $\ell\in L$, for each $j\in J$.
The hybrid map is called \emph{constant-rank}, \emph{injective}, or \emph{surjective} if each $f|_{M_j}$ is as well.
It is called an \emph{embedding} if each $f|_{M_j}$ is an embedding and $f$ is a homeomorphism onto its image.
}

\defn{\label{def:hpf}
The \emph{hybrid pushforward} $f_*:TM\into TN$ is the hybrid map defined piecewise as $f_*|_{TM_j} := (f|_{M_j})_*$.
}


\defn{\label{def:hvf}
A \emph{hybrid vector field} on a hybrid manifold $M:=\coprod_{j\in J}M_j$ is a hybrid map $G:M\into TM$ for which $G|_{M_j}$ is a vector field on $M_j$, i.e. $G|_{M_j} \in \e{T}(M_j)$.
We let $\e{T}(M)$ denote the space of hybrid vector fields on $M$.
}

To state the main result of this paper, we need to embed manifolds into hybrid manifolds.  
This can be achieved by first partitioning the smooth manifold to obtain a hybrid manifold, then embedding this hybrid manifold via the previous definitions.

\defn{\label{def:part} 
A \emph{partition} of an $n$-dimensional manifold $M$ is a finite set $\set{M_j}_{j\in J}$ of embedded $n$-dimensional submanifolds $M_j\subset M$ for which $\bigcup_{j\in J} M_j = M$ and if $i\ne j$ we have $\Int(M_j)\cap \Int(M_i) = \emptyset$.}

\defn{\label{def:hembed}
A \emph{hybrid embedding} of a manifold $M$ into a hybrid manifold $N := \coprod_{j\in J} N_j$ is determined by a partition $\set{M_j}_{j\in J}$ of $M$ and a hybrid embedding 
\eqn{f:\coprod_{j\in J}M_j\into \coprod_{j\in J} N_j}
for which $f_j : M_j \into N_j$ for each $j\in J$.  Any $G\in\e{T}(M)$ may be pushed forward to a unique $f_* G\in\e{T}(f(M))$.
The image of $f$ is a \emph{hybrid embedded submanifold}.
}

With these preliminaries established, we can define the class of hybrid systems considered in this paper.
To the best of our knowledge, this definition of a hybrid dynamical system has not appeared before.
However, in light of the constructions contained in this section,
it may be seen as a mild generalization of a \emph{simple hybrid system} \cf{\S 3.2 in~\cite{Ames2006}}.  
Further, in Section~\ref{sec:inv} we will see that this definition supports powerful geometric analysis of the dynamics near a hybrid periodic orbit.
Finally, this class of hybrid systems encompasses many closed-loop models of legged locomotion~\cite{McGeer1990, GrizzleAbba2002, CollinsRuina2005, GhigliazzaAltendorfer2003, ProctorHolmes2008, WesterveltGrizzle2003, MorrisGrizzle2009, PoulakakisGrizzle2009, AmesGregg2006}.
We contend that these facts justify the introduction of the novel definition.

\defn{A \emph{hybrid dynamical system} is specified by a triple $H := (D,F,R)$ where:
\begin{itemize}
\item[$D$] $= \coprod_{j\in J} D_j$ is a hybrid manifold;
\item[$F$] $\in\e{T}(D)$ is a hybrid vector field on $D$;
\item[$R$] $:S\into T$ is a hybrid map, $S,T\subset\bd D$ are hybrid embedded submanifolds, and $S$ has codimension~1.
\end{itemize}
As in~\cite{Ames2006}, we call $R$ the \emph{reset map} and $S$ the \emph{guard}.
}

Note that if $F$ is tangent to $S$ at $x\in D$, there is a possible ambiguity in determining a trajectory from $x$---one may either follow the flow of $F$ on $D$ or apply the reset map to obtain a new initial condition $y = R(x)$.

\assump{\label{assump:transS}
To ensure that trajectories are uniquely defined, we assume that $F$ is outward-pointing on $S$.
}

\rem{As defined above, hybrid dynamical systems possess unique \emph{executions} or \emph{trajectories} from every initial condition.  
This fact can be demonstrated algorithmically.
For any $x\in D_j$, obtain the maximal integral curve of $F|_{D_j}$.
This integral curve must either: a) continue for all time; b) exit $D_j$ without intersecting $\bd D_j$ (in which case execution terminates); or c) intersect the boundary at $y\in\bd D_j$.  
If $y\in S$, the map $R$ is applied to obtain a new initial condition $R(y)\in T$, and otherwise execution terminates.}

The following definition enables us to embed smooth dynamical systems into hybrid dynamical systems in such a way that trajectories of the smooth system are preserved in the hybrid system.
We illustrate the use of this construction by giving a terse description of trajectories for this class of hybrid systems.
In the subsequent sections, we use this construction to state the main results of this paper.

\defn{\label{def:dsembed}
A \emph{hybrid dynamical embedding} of a dynamical system $(M,G)$ into a hybrid dynamical system $(D,F,R)$ is a hybrid embedding $f:M\into D$ for which $f_* G = F|_{f(M)}$ and $R|_{f(M)\cap S}$ is a hybrid diffeomorphism from $f(M)\cap S$ onto $f(M)\cap T$.}

\rem{A trajectory of a hybrid dynamical system H may be obtained from a hybrid dynamical embedding of the system $\paren{J,\vf{t}}$, where $J\subset\R$ is a connected interval.}

\defn{A \emph{$\tau$-periodic orbit} of a hybrid dynamical system is a hybrid dynamical embedding $\gamma$ of the dynamical system $(S^1,\frac{2\pi}{\tau} \vf{\vphi})$, where $S^1$ is the unit circle.}

\rem{We alternately refer to $\gamma$ as a \emph{periodic trajectory} and often write $\gamma$ in place of the image $\gamma(S^1)$.}

\subsection{Hybrid \Pmap}\label{sec:P}

To state the main result of this paper, we must construct the Poincar\'{e} map associated with a periodic orbit of a hybrid system.  
This has been developed before~\cite{MorrisGrizzle2009, WendelAmes2010}; the construction is more delicate than for smooth systems since trajectories of hybrid systems do not necessarily vary continuously with initial conditions.  
We directly demonstrate this continuous dependence in the construction of the map.

Let $H = (D,F,R)$ be a hybrid dynamical system and $\gamma$ a periodic orbit of $H$ with period $\tau$.  
Then $\gamma$ undergoes a finite number of transitions $k\in\N$, so we may index the corresponding sequence of domains as\footnote{We regard subscripts modulo $k$ so that $D_k \equiv D_0$.} $D_1,\dots,D_k$.
Without loss of generality, assume the $D_j$'s are distinct\footnote{Otherwise we can find $\set{B_j}_{j=1}^k$ be such that $B_j\subset D_j$ is open, $\gamma\subset\bigcup_{j=1}^k B_j$, and $B_i\cap B_j=\emptyset$ if $i\ne j$, then proceed on $D := \coprod_{j=1}^k B_j$.}; let $\gamma_j := \gamma\cap T_j$ be the entry point of $\gamma$ in $D_j$ and let $\tau_j$ be the time spent by $\gamma$ in $D_j$.
We wish to construct the Poincar\'{e} map $P$ associated with $\gamma$ over a neighborhood of $\gamma_j$ in $T$.
To do this, we must ensure that each initial condition in that neighborhood has a well-defined non-zero first-return time to $T$; the following assumption guarantees this.
\assump{\label{assump:transT}
To ensure the Poincar\'{e} map is well-defined, we assume $F$ is transverse to $T$ and not outward-pointing.
}

\begin{figure}[ht]
\centering
\includegraphics[scale=0.80]{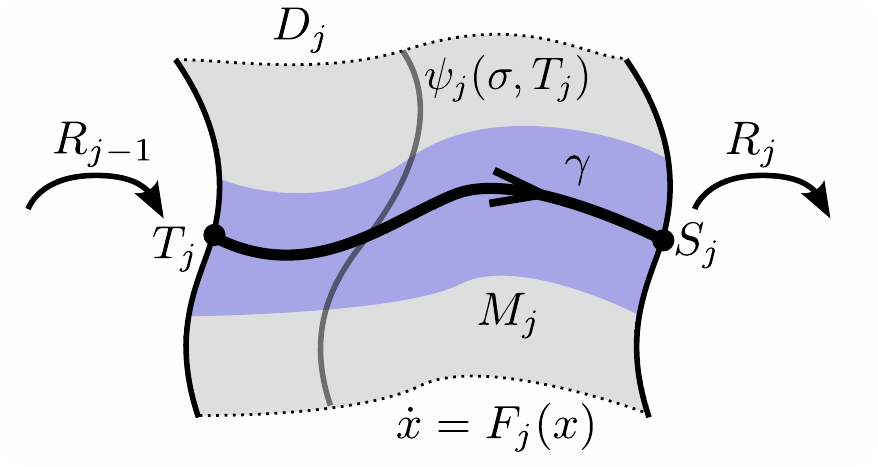}
\caption{
Illustration of a hybrid domain containing part of a periodic orbit.
The periodic orbit $\gamma$ enters the domain $D_j$ inside the submanifold $T_j\subset\bd D_j$ via the reset map $R_{j-1}$.
Initial conditions in $T_j$ flow to the submanifold $S_j\subset\bd D_j$ via the vector field $F_j$.
The map $\psi_j$ smoothly reparameterizes these trajectories by the time required to pass from $T_j$ to $S_j$, thus $\psi_j(\sigma,T_j)$ is an embedded submanifold for all $\sigma\in[0,1]$ and $S_j = \psi_j(1,T_j)$.
While in domain $D_j$, $\gamma$ lies in the invariant submanifold $M_j$ constructed in Theorem~\ref{thm:inv}.
By construction, $M_j$ is an integral submanifold of $F_j$
and $\dim M_j \le \dim D_j$; see Fig.~\ref{fig:result} for an illustration when $\dim M_j< \dim D_j$.
}
\label{fig:domain}
\end{figure}

\noindent Now for $j = 1,\dots,k$ and referring to Fig.~\ref{fig:domain} for an illustration of these objects, let:

\begin{enumerate}
    \item[$\phi_j$] $:\e{F}_j\into D_j$ be the maximal flow of $F$ on $D_j$;
    \item[$T_j$] $\subset T\cap D_j$ be a neighborhood of $\gamma_j$ over which Lemma~\ref{lem:sec} may be applied between $T$ and $S$ on $D_j$;
    \item[$\psi_j$] $:[0,1]\times T_j\into D_j$ be the embedding from Lemma~\ref{lem:sec};
    \item[$S_j$] $:= \psi_j(1,T_j)\subset S\cap D_j$ be the image of $T_j$ in $S$ under the flow on $D_j$;
    \item[$R_j$] $:S_j\into T$ denote the restriction $R_j := R|_{S_j}$;
    \item[$p_j$] $: T_j\into T$ be defined by $p_j(u) := R_j(\psi_j(1,u))$.
\end{enumerate}
The Poincar\'{e} map over the section $T_j$ is obtained formally by iterating the $p$'s around the cycle:  
\eqnn{\label{eq:P} P_j := p_{j-1}\circ\cdots\circ p_1\circ p_k\circ\cdots\circ p_j.  }
The neighborhood $\Sigma_j\subset T_j$ of $\gamma_j$ over which this map is well-defined is determined by pulling $T_j$ backward around the cycle,
\eqn{\Sigma_j = \paren{p_j^{-1}\circ\cdots\circ p_k^{-1}\circ p_1^{-1}\circ\cdots\circ p_{j-1}^{-1}}(T_j),}
and similarly for any iterate of $P_j$.  

It is a standard result for smooth dynamical systems that \emph{Floquet multipliers} (the eigenvalues of the linearized \Pmap) do not depend on the choice of Poincar\'{e} section \cf{Section 1.5 in~\cite{GuckenheimerHolmes1983}}.
The following lemma generalizes this result to the hybrid setting by demonstrating that if a \Pmap\, obtained from one domain has an attracting invariant submanifold via Lemma~\ref{lem:nondeg}, then the map obtained in any other starting domain has a diffeomorphic attracting submanifold.
As a consequence, non-zero Floquet multipliers are shared between the $P_j$'s after a sufficient number of iterations.

\lem{\label{lem:pmap}
Let $j\in\set{1,\dots,k}$ and $n \ge \min_\ell\dim D_\ell$.  If $P_j^n$ has constant rank equal to $r$ near $\gamma_j$, 
then $P_\ell^{n+1}$ has constant rank equal to $r$ near $\gamma_\ell$ for all $\ell\in\set{1,\dots,k}$.
}

In the proof of Lemma~\ref{lem:pmap}, we make use of an elementary fact from linear algebra.
The result is easily obtained from Sylvester's inequality~\cf{Appendix A.5.4 in~\cite{CallierDesoer1991}}.
\prop{\label{prop:la2}
For $j\in\set{1,\dots,k}$, 
suppose $A_j\in\R^{n_{j+1}\times n_j}$ where $n_k = n_1$, 
define $B_j := A_{j-1}\cdots A_1 A_k \cdots A_j$, 
and let $n \ge \min_j n_j$.
Then for all $\ell\in\set{1,\dots,k}$, we have $\rk B_\ell^{n+1} = \rk B_j^{n+1} = \rk B_j^n$.
}

\pf{{\it (of Lemma~\ref{lem:pmap})}
By Lemma~\ref{lem:nondeg}, there is a neighborhood $N_j\subset T_j$ of $\gamma_j$ on which $P_j^n$ has constant rank equal to $r$.
Fix $\ell\in\set{1,\dots,k}$, let $p_\ell^j~:=~p_{j-1}~\circ~\cdots~\circ~p_\ell$, and define $N_\ell := (p_\ell^j)^{-1}(N_j)$.
Then $N_\ell\subset T_\ell$ is a neighborhood of $\gamma_\ell$ and, for all $x\in N_\ell$, 
by Sylvester's inequality 
\[ \rk (P_\ell^{n+1})_*(x) \le \rk (P_j^n)_*\paren{p_\ell^j(x)} = r. \]
Furthermore by Proposition~\ref{prop:la2},
\eqnal{
  \rk (P_\ell^{n+1})_*(\gamma_\ell) & = \rk (P_j^{n+1})_*(\gamma_j) \\
  & = \rk (P_j^n)_*(\gamma_j) = r.
}
We conclude the rank of $P_\ell^{n+1}$ is at least $r$ on a neighborhood $L_\ell\subset T_\ell$ of $\gamma_\ell$,
whence $\rank P_\ell^{n+1} = r$ on $L_\ell\cap N_\ell$.
}
As a consequence, if the \Pmap\, associated with \emph{any} section for the periodic orbit $\gamma$ satisfies the hypotheses of Lemma~\ref{lem:nondeg}, then the \Pmap\, associated with \emph{any other} section also satisfies the hypotheses.

\rem{\label{rem:pmap}
It may be easier to evaluate the rank of the \Pmap\, in some domains than others.  In particular, if $P_j$ is a diffeomorphism for some $j\in\set{1,\dots,k}$, then all iterates are constant rank.
}

\subsection{Hybrid Invariant Subsystem}\label{sec:inv}

This section contains the main result of this paper:  
when iterates of the \Pmap\, associated with a periodic orbit of a hybrid dynamical system have constant rank,
trajectories starting near the orbit converge in finite time to an embedded smooth dynamical system.

\thm{\label{thm:inv}
Let $H = (D,F,R)$ be a hybrid dynamical system, $\gamma$ a periodic orbit of $H$, and suppose the composition of any \Pmap\, for $\gamma$ with itself at least $\min_j\dim D_j$ times has constant rank equal to $r$ on a neighborhood of its fixed point.  
Then there is an $(r+1)$-dimensional dynamical system $(M, G)$, a hybrid dynamical embedding $f:M\into D$, and an open hybrid set $W\subset D$ so that
$\gamma\subset f(M)\cap W$ and
trajectories starting in $W$ flow into $M$ in finite time.
}

\pf{
By assumption, we may apply Lemma~\ref{lem:nondeg} to $P$ to obtain a neighborhood $N\subset T_1$ of $\gamma\cap T_1$, an embedded submanifold $\Sigma\subset T_1$ containing $\gamma\cap T_1$, and a pair of neighborhoods $U,V\subset\Sigma$ of $\gamma\cap T_1$ so that $P|_U:U\into V$ is a diffeomorphism and $\dim U = r$.  
Now we consider the subset of $D$ obtained by propagating each $x\in U$ around one cycle.
Let $U_1 = U$ and $U_j = P_{j-1}(U_{j-1})$ for $j = 2,\dots, k$. 
Away from the boundaries, we can obtain the desired set directly from $\psi_j$ as $\Int(M_j) := \psi_j( (0,1),U_j)$.
For $j = 2,\dots,k-1$ we can simply attach the corresponding boundaries to obtain $M_j := \psi_j([0,1],U_j)$.
However, since we may not assume $U\subset V$ or $V\subset U$ (only that $U\cap V$ is a neighborhood of $\gamma\cap T_1$), we must be careful in attaching the boundary between $M_k$ and $M_1$.
Thus, we let $M_1 = \psi_1((0,1],U_1)\cup (U_1\cap P_k(U_k))$ and $M_k = \psi_k([0,1),U_k)\cup (\psi_k(1,U_k)\cap R_k^{-1}(U_1))$.
With this construction, for each $j=1,\dots,k$ we have that $M_j$ is a smooth submanifold with boundary $\bd M_j\subset T_j\cup S_j$ and $\bd M_j$ contains both points in $\gamma_j\cap\bd D_j$; see Fig.~\ref{fig:domain} for an illustration of $M_j$.

Since $M_j$ is an integral submanifold of $F$ on $D_j$, the vector field $F$ restricts to $M_j$.
Letting $G_j$ denote this restriction, each $(M_j,G_j)$ is a smooth dynamical system and $G_j$ points inward on $\bd M_j\cap T_j$ and outward on $\bd M_j\cap S_j$.  
Since $P|_U$ is a diffeomorphism, each $R_j|_{\bd M_j\cap S_j}:\bd M_j\cap S_j\into \bd M_{j+1}\cap T_{j+1}$ is a diffeomorphism as well.  
Therefore we may glue these systems together one-by-one via Lemma~\ref{lem:glue} to obtain a smooth dynamical system without boundary $(M,G)$ which embeds into $H$ and contains $\gamma$.

Finally, let $N,U,$ and $V$ be as above and let $\delta > 0$ be an arbitrary positive number.  
Note that by continuity of $P$ and the time-to-impact maps of Lemma~\ref{lem:sec}, there is a neighborhood $W_1\subset N\subset T_1$ of $\gamma\cap T_1$ so that $P^n(W_1)\subset U\cap V$ and each $w\in W_1$ flows into $U\cap V$ before time $n\tau+\delta$.  
Since the $P_j$'s are continuous, for $j=2,\dots,k$ there are neighborhoods $W_j\subset T_j$ of $\gamma\cap T_j$ so that every $w\in W_j$ flows into $W_{j+1}$ before time $\tau_j+\delta/k$.
Taking the union of these neighborhoods as $W = \coprod_{j=1}^k\psi_j([0,1],W_j)$ yields an open hybrid submanifold $W\subset D$ so that $\gamma\subset M\cap W$ and every point in $W$ flows into $W_1$ before time $\tau+\delta$, and hence into $M$ before time $(n+1)\tau + 2\delta$; see Fig.~\ref{fig:result} for an illustration of these neighborhoods in a particular two-domain hybrid dynamical system.
}

\cor{\label{cor:hstab}
$\gamma$ is asymptotically stable for $H = (D,F,R)$ if and only if $\gamma$ is asymptotically stable for $(M,G)$.
}

\pf{
Since all trajectories in a neighborhood $W$ of $\gamma$ reach $M$ in finite time and the hybrid flow is continuous near $\gamma$, trajectories in $W$ will converge to $\gamma$ asymptotically if and only if trajectories in $W\cap M$ converge to $\gamma$ asymptotically.
This occurs precisely when $\gamma$ is asymptotically stable for $(M,G)$ since by construction $M$ is an integral submanifold of $F$ and $R|_{M\cap S} : M\cap S\into M\cap T$ is a diffeomorphism. 
}

If each of the $D_j$'s have the same dimension and $R:S\into T$ is a diffeomorphism, the rank condition of Theorem~\ref{thm:inv} is trivially satisfied, and we can globalize the construction using Lemma~\ref{lem:glue}.  
This provides a smooth $n$-dimensional generalization of the construction in~\cite{SimicJohansson2005}. 

\cor{\label{cor:diffeo}
Let $H = (D,F,R)$ be a hybrid dynamical system with $D = \coprod_{j\in J}D_j$, $R:S\into T$, and $\bd D = S\cup T$.  If $\dim D_j = n$ for all $j\in J$ and $R$ is a diffeomorphism, then there is a surjective hybrid dynamical embedding from an $n$-dimensional dynamical system $(M,G)$ onto $H$.
}

\section{Examples}
\label{sec:ex}

\subsection{Hybrid Floquet Coordinates}

The following single-domain system clearly satisfies the hypotheses of Theorem~\ref{thm:inv}, and demonstrates the canonical form for hybrid Floquet coordinates.

\ex{
Let $H = (D,F,R)$ be a hybrid system over the single domain $D = [0,1]\times\R^k\times\R^\ell$ with 
vector field $F(t,x,z) = \vf{t} + \sum_{j=1}^k f^j(t,x)\vf{x^j} + \sum_{i=1}^\ell g^i(t,z^i) \vf{z^i}$, 
reset map $R:\set{1}\times\R^k\times\R^\ell\into\set{0}\times\R^k\times\R^\ell$ defined by $R(1,x,z) = (0,x,Az)$ where $A\in\R^{\ell\times\ell}$ is nilpotent, $f(t,\xi) = 0$ for all $t\in[0,1]$ and some $\xi\in\R^k$, and $g^i(t,0) = 0$ for all $i$.  
Consider the Poincar\'{e} map 
\eqn{P:\set{0}\times\R^k\times\R^\ell\into\set{0}\times\R^k\times\R^\ell.}
It is clear that $P(0,\xi,0) = (0,\xi,0)$,
\eqn{P^\ell(\set{0}\times\R^k\times\R^\ell) = \set{0}\times\R^k\times\set{0},}
\eqn{\rk R|_{\set{1}\times\R^k\times\set{0}} = k.}
Therefore $\rk P^{k+\ell} = k$, whence we may apply Theorem~\ref{thm:inv}.
The resulting smooth invariant subsystem is diffeomorphic to $S^1\times\R^k$.
}

\begin{figure}
  \centering
  \subfloat[vertical hopper schematic]{\label{fig:hopsch}
    \begin{minipage}[l]{0.45\columnwidth}
    \includegraphics[height=4.5cm]{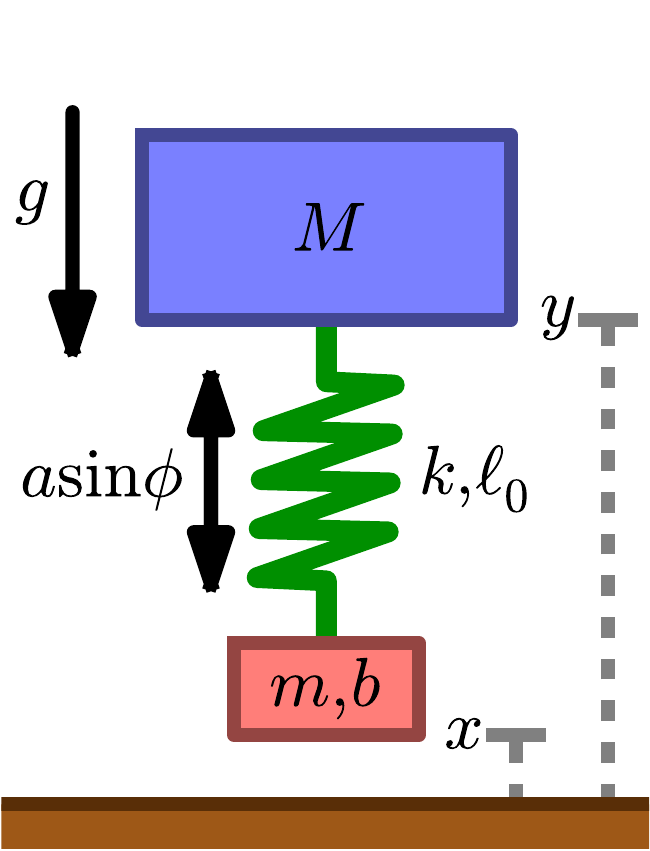}
    \end{minipage}
  }
  \hspace{0.2cm}
  \subfloat[hopper dynamics]{\label{fig:hopeqs}
    \begin{minipage}[r]{0.4\columnwidth}
    \footnotesize
    For $(\phi,x,\dot{x},y,\dot{y}) \in D_a$:
    \eqnal{
      \dot{\phi} =\, & \omega, \\
      m\ddot{x} =\, & - k\ell_0 - a\sin\phi \\
      & + k (y-x) - b\dot{x} - g m, \\
      M\ddot{y} =\, & k\ell_0 + a\sin\phi \\
      & - k (y-x) - g M. \\
    }
    For $(\phi,y,\dot{y}) \in D_g$:
    \eqnal{
      \dot{\phi} =\, & \omega, \\
      M\ddot{y} =\, & k\ell_0 + a\sin\phi \\
      & - k y - g M.
    }
    \end{minipage}
  }
  \caption{
  Schematic and dynamics of vertical hopper.
  (a) Two masses $m$ and $M$, constrained to move vertically above a ground plane in a gravitational field with strength $g$, are connected by a linear spring with stiffness $k$ and nominal length $\ell_0$.
The spring is equipped with an actuator that exerts a phase-varying force $a\sin\phi$ where $\dot{\phi} = \omega$.
  The lower mass experiences viscous drag proportional to velocity with constant $b$ when it is in the air, and impacts plastically with the ground (i.e. it is not permitted to penetrate the ground and its velocity is instantaneously set to zero whenever a collision occurs).
  (b) Hopper dynamics in aerial and ground domains.
  }
\label{fig:hop}
\end{figure}

\subsection{Vertical Hopper}

We apply Theorem~\ref{thm:inv} to demonstrate the existence of low-dimensional invariant dynamics in the model for forced vertical hopping illustrated in Fig.~\ref{fig:hopsch}.
The state space in the aerial phase is $D_a := S^1 \times T\R_{\ge 0}\times T\R$.
Writing $(\phi,x,\dot{x},y,\dot{y})\in D_a$, 
the aerial dynamics are given in Fig.~\ref{fig:hopeqs}.
When the lower mass rests on the ground, the state space resides in $D_g := S^1\times T\R$ and the dynamics of the upper mass are obtained by restricting to the submanifold $\set{(\phi,x,\dot{x},y,\dot{y}) : x = \dot{x} = 0}$ as in Fig.~\ref{fig:hopeqs}.
Transition from the aerial to the ground domain occurs when the lower mass collides with the ground, and the state is reset according to $(\phi,0,\dot{x},y,\dot{y}) \mapsto (\phi,y,\dot{y})$.
The lower mass lifts off when the normal force required to keep it from penetrating the ground plane becomes zero, i.e. when $m g = -k\,\ell_0 - a\sin\phi + k\,y$, and the state is reset via $(\phi,y,\dot{y}) \mapsto (\phi,0,0,y,\dot{y})$.

Numerical simulations\footnote{
Note that simulation of hybrid dynamical systems is non-trivial.
We make use of a recently-developed algorithm with desirable convergence properties~\cite{BurdenGonzalez2011}.
In particular, we use Euler step size $h = 1\times 10^{-3}$ and relaxation parameter $\eps = 1\times 10^{-12}$.
As a note to practitioners, we found that numerical linearization of the \Pmap\, via finite differenes was sensitive to the coordinate displacement when using large values for the relaxation parameter.
The sourcecode for this simulation is available online at {\tt http://purl.org/sburden/cdc2011}
}
 indicate that with parameters $(m,M,k,b,\ell_0,a,\omega,g) = (1,2,10,5,2,20,2\pi,2)$, the hybrid system possesses a stable periodic orbit, $\gamma$.
Choosing a Poincar\'{e} section in domain $D_g$ at $\phi = 0$, we find that $\gamma$ intersects this section at the point $(y,\dot{y}) = (1.96, 1.88)$ and that the eigenvalues of the linearized \Pmap\, are $-0.25\pm0.70j$.
Both eigenvalues lie inside the unit disc, corroborating the observed stability of the orbit.
Further, since neither eigenvalue is close to zero, we conclude the \Pmap\, has full rank equal to 2 near its fixed point.
Therefore by Remark~\ref{rem:pmap} the hypotheses of Theorem~\ref{thm:inv} are satisfied, and we conclude the system's dynamics collapse to a smooth 3-dimensional subsystem after one hop.

\addtolength{\textheight}{-0.0cm}

\section{Discussion}\label{sec:disc}
We demonstrated the existence of a locally attracting constant-dimensional invariant subsystem near a hybrid periodic orbit whenever iterates of the associated \Pmap\, have constant rank.  
Under a genericity condition, near a periodic orbit of a smooth dynamical system there exist \emph{Floquet coordinates} in which the dynamics decouple into a constant-frequency phase variable and a time-invariant transverse linear system~\cite{floquet1883edl, Guckenheimer1975, GuckenheimerHolmes1983}.
Under the additional rank hypothesis of Theorem~\ref{thm:inv}, we obtain a canonical form for the Floquet structure of a hybrid periodic orbit.
Indeed, the smooth subsystem $(M,G)$ foliates the dynamics near the periodic orbit in each domain.  
Thus the behavior of the hybrid system near the orbit is a trivial extension of the behavior of a smooth system---portions of the smooth dynamics are ``stacked'' in transverse coordinates and annihilated within a finite number of cycles via a nilpotent linear operator. 
On the smooth subsystem, the standard construction of Floquet coordinates may be applied, 
generalizing the class of systems which may be analyzed using the empirical approach developed in~\cite{Revzen2009, RevzenGuckenheimer2011}.

In addition to providing a canonical form for the dynamics near such non-degenerate periodic orbits, the results of this paper suggest a mechanism by which a many-legged locomotor may formally collapse a large number of degrees-of-freedom to produce a low-dimensional coordinated gait.
This provides a link between currently disparate lines of research, namely the formal analysis of hybrid periodic orbits, the design of robots for locomotion and manipulation tasks, and the scientific probing of neuromechanical control architectures in organisms.
It shows that hybrid models naturally exhibit dimension reduction, that this reduction may be deliberately designed into an engineered system, and that evolution may have exploited this reduction in developing its spectacular locomotors.

%
%


\small{

\subsection*{Acknowledgements \& Support}
We thank Saurabh Amin, Jonathan Glidden, Humberto Gonzalez, John Guckenheimer, and Ramanarayan Vasudevan for helpful conversations and careful readings of this paper.

S. Burden was supported in part by an NSF Graduate Research Fellowship.
S. Revzen was supported in part by NSF Frontiers for Integrative Biology Research (FIBR), Grant No. 0425878-Neuromechanical Systems Biology.
Part of this research was sponsored by the Army Research Laboratory under Cooperative Agreements W911NF-08-2-0004 and W911NF-10-2-0016.
}

\bibliographystyle{unsrt}
\bibliography{hsfloq}




\end{document}